\documentclass{amsart}
\usepackage{amssymb}
\usepackage{amsmath}
\usepackage{amsthm}
\usepackage{color,epsfig}
\usepackage{graphicx}
\usepackage{mathdots}

\usepackage[dvipdfm,colorlinks,backref,pagebackref,hypertexnames=false,linkcolor=blue,urlcolor=green,citecolor=red]{hyperref}
\newcommand{\R}{{\mathbb{R}}}
\newcommand{\C}{{\mathbb{C}}}
\newcommand{\F}{{\mathbb{F}}}

\newcommand{\la}{\langle}
\newcommand{\ra}{\rangle}

\newcommand{\al}{\alpha}
\newcommand{\lam}{\lambda}
\newcommand{\Sig}{\Sigma}

\newcommand{\Irr}{{\rm Irr}}
\newcommand{\GL}{{\rm GL}}

\renewcommand{\epsilon}{\varepsilon}

\newtheorem{theorem}{Theorem}[section]
\newtheorem{lemma}[theorem]{Lemma}
\newtheorem{proposition}[theorem]{Proposition}
\newtheorem{corollary}[theorem]{Corollary}
\theoremstyle{remark}
\newtheorem{definition}[theorem]{Definition}
\newtheorem{example}[theorem]{Example}
\newtheorem*{remark}{Remark}

\begin{document}

\title[Supercharacters and pattern subgroups]{Supercharacters and pattern subgroups in the upper triangular groups}

\author{Tung Le}

%\address{Institute of Mathematics, University of Aberdeen, Aberdeen AB24 3UE, Scotland, UK}

%\email{t.le@abdn.ac.uk}
\email{lttung96@yahoo.com}

\date{March 2nd, 2010}



%% 20C15   	Ordinary representations and characters
%% 20D15   	Nilpotent groups, $p$-groups
%% 20D20   	Sylow subgroups, Sylow properties, $\pi$-groups, $\pi$-structure
%% 20C33   	Representations of finite groups of Lie type

\maketitle

%%%%%%%%%%%%%%%%%%%%%%%%%%%%%%%%%%%%%%%%%%%%%%%%%%%%%%%%%%%%%%%%%%%%%%%%%%%%%%%%%%%%%%%%%%

\begin{abstract}

Let $U_n(q)$ denote the upper triangular group of degree $n$
over the finite field $\F_q$ with $q$ elements. It is known that irreducible constituents of supercharacters partition  the set of all irreducible characters $\Irr(U_n(q))$. In this paper we present a correspondence between supercharacters and pattern subgroups of the form $U_k(q)\cap {}^wU_k(q)$ where $w$ is a monomial matrix in $\GL_k(q)$ for some $k<n$.

\end{abstract}

%%%%%%%%%%%%%%%%%%%%%%%%%%%%%%%%%%%%%%%%%%%%%%%%%%%%%%%%%%%%%%%%%%%%%%%%%%%%%%%%%%%%%%%%%%

\section{Introduction}
\label{sec:intro}
Let $q$ be a power of a prime $p$ and $\F_q$ a field with $q$
elements. The group $U_n(q)$ of all upper triangular $(n\times
n)$-matrices over $\F_q$ with all diagonal entries equal to $1$ is a
Sylow $p$-subgroup of $\GL_n(\F_q)$. It is conjectured by  Higman
\cite{Hig} that the number of conjugacy classes of $U_n(q)$ is given
by a polynomial in $q$ with integer coefficients. Isaacs \cite{Is2} shows that the degrees of all
irreducible characters of $U_n(q)$ are powers of $q$. Huppert \cite{Hup} proves that character degrees of $U_n(q)$ are precisely of the form $\{q^e:0\leq e \leq \mu(n)\}$ where the upper bound $\mu(n)$ was known to Lehrer \cite{Leh}.
Lehrer \cite{Leh} conjectures that each number  $N_{n,e}(q)$ of
irreducible characters of $U_n(q)$ of degree $q^e$ is given by a
polynomial in $q$ with integer coefficients. Isaacs \cite{Is3}
suggests a strengthened form of Lehrer's conjecture stating that
$N_{n,e}(q)$ is given by a polynomial in $(q-1)$ with nonnegative
integer coefficients. So, Isaacs' Conjecture implies Higman's and
Lehrer's Conjectures.

Many efforts have been made to understand more about $U_n(q)$, see \cite{Carlos2,Vera-Lopez,Is-Di,Anton1,Is2,Is3,Robinson,Thompson1,Thompson2}, among others.
%%Thompson \cite{Thompson1}, Robinson \cite{Robinson}, Andr\'e \cite{Carlos2}, Isaacs \cite{Is2,Is3}, Diaconis-Isaacs \cite{Is-Di}, Arregi-Vera-L\'opez \cite{Vera-Lopez}, Evseev \cite{Anton1}...
Supercharacters arise as tensor products of some elementary characters to give a `nice' partition of all non-principal irreducible characters of $U_n(q)$, see \cite{Carlos2,Tung1}. Supercharacters have been defined for Sylow $p$-subgroups of other finite groups of Lie type (see \cite{Carlo-Neto}), and in general for finite algebra groups (see \cite{Is-Di}).

Here, for $U_n(q)$ we show a natural correspondence between supercharacters and pattern subgroups (Theorem \ref{main_thm}). To highlight the main idea of construction, we have deferred all of our proofs to Section \ref{sec:proofs}. Next we apply this correspondence to decompose certain supercharacters into irreducible constituents in the last section.

\section{Suppercharacters and pattern subgroups}
\label{sec:spchr+ptgp}

Let $\Sig:=\Sig_{n-1}=\la \alpha_1,...,\alpha_{n-1}\ra$ be the root system
of $\GL_n(q)$ with respect to the maximal split torus equal to the
diagonal group, see \cite[Chapter 3]{cart2}. Denote $\alpha_{i,j}:=\alpha_i+\alpha_{i+1}+...+\alpha_j$ for all $0<i\leq j<n$. Denote by $\Sig^+$ the set of all positive roots.  The root subgroup $X_{\al_{i,j}}$ is the set of all matrices of the form $I_n+c\cdot e_{i,j+1}$, where $I_n=$ the  identity $(n\times n)$-matrix, $c\in\F_q$ and $e_{i,j+1}=$ the zero matrix except a `$1$' at entry $(i,j+1)$. The upper triangular group $U_n(q)$ is generated by all $X_\al$ where $\al\in \Sig^+$. We write $U$ for $U_n(q)$ if $n$ and $q$ are clear from the context. For convenience when using the root system, we consider the upper triangular group as a tableaux.
\begin{small}
\[
\left(\begin{array}{ccccc }
1 & * & * & * & *    \\
. & 1 & * & * & *   \\
. & . & 1 & * & *   \\
. & . & . & 1 & *    \\
. & . & . & . & 1
\end{array}\right)
\longrightarrow
\begin{array}{ccccc}
\cline{2-5}
& \multicolumn{1}{|c|}{\alpha_1}&\multicolumn{1}{|c|}{\alpha_{1,2}}& \multicolumn{1}{|c|}{\alpha_{1,3}}& \multicolumn{1}{|c|}{\alpha_{1,4}}  \\
\cline{2-5}
  &  & \multicolumn{1}{|c|}{\alpha_2} &\multicolumn{1}{|c|}{\alpha_{2,3}} & \multicolumn{1}{|c|}{\alpha_{2,4}}   \\
\cline{3-5}
 &   &  & \multicolumn{1}{|c|}{\alpha_3} & \multicolumn{1}{|c|}{\alpha_{3,4}}  \\
 \cline{4-5}
  &   &   &   & \multicolumn{1}{|c|}{\alpha_4} \\
  \cline{5-5}
\end{array}
\]
\end{small}

A subset $S\subset \Sig^+$ is called {\em closed} if for each $\al,\beta\in S$ such that $\al+\beta\in\Sig^+$ then $\al+\beta\in S$. A {\em pattern} subgroup of $U$ is a group generated by all root subgroups $X_\al$, where $\al\in S$ a closed positive root subset.

Let $G$ be a group. Denote $G^\times:=G-\{1\}$, $\Irr(G)$ the set of all complex irreducible characters of
$G$, and $\Irr(G)^\times:=\Irr(G)-\{1_G\}$. For $H\unlhd G$, let $\Irr(G/H)$ denote the set of all irreducible characters of $G$ with $H$ in the kernel. If $K\leq G$ such that $G=H\rtimes K$, then for each character $\xi$ of $K$, we denote the inflation of $\xi$ to $G$ by $\xi_G$, i.e. $\xi_G$ is the extension of $\xi$ to $G$  with $H\subset \ker(\xi_G)$.
Furthermore, for $H\leq G$ and $\xi\in \Irr(H)$, we define $\Irr(G,\xi):=\{\chi\in \Irr(G):(\chi,\xi^G)\neq 0\}$ the irreducible constituent set of $\xi^G$, and for $\chi\in \Irr(G)$, we denote its restriction to $H$ by $\chi|_H$.

For a field $K$, let $K^\times:=K-\{0\}$. In the whole paper, we fix a
nontrivial linear character $\varphi:\F_q\to \C^\times$. For each $\al\in \Sigma^+$ and $s \in \F_q$, the
map $\phi_{\al,s}: X_\al\to \C^\times$, $x_\al(d)\mapsto \varphi(ds)$ is a linear character of $X_\al$,
and all linear characters of $X_\al$ arise in this way.

For each $\al_{i,j}$, we define
\[arm(\al_{i,j}):=\{\al_{i,k}:i\leq k < j\}  \mbox{ and } leg(\al_{i,j}):=\{\al_{k,j}:i< k\leq j\}.\]
If $i=j$,  $\al_{i,i}=\al_i$, then $arm(\al_i)$ and $leg(\al_i)$ are empty. For each $\al\in \Sig^+$, we define the {\em hook} of $\al$ as $h(\al):=arm(\al)\cup leg(\al)\cup\{\al\}$, the {\em hook group} of $\al$ as $H_\al:=\langle X_\beta:\beta\in h(\al)\rangle$, and the {\em base group} $V_\al:=\langle X_\beta:\beta\in \Sig^+-arm(\al)\rangle$. Since $[V_\al,V_\al]\cap X_\al=\{1\}$, for each $s\in\F_q^\times$ there exists a linear $\lam_{\al,s}\in \Irr(V_\al)$ such that $\lam_{\al,s}|_{X_\al}=\phi_{\al,s}$ and $\lam_{\al,s}|_{X_\beta}=1_{X_\beta}$ for the others $X_\beta\subset V_\al$, $\beta\neq \al$. Denote by $\Irr(V_\al/[V_\al,V_\al])^\times$ the set of all these linear characters of $V_\al$.

\begin{lemma}
$\lam_{\al,s}^U$ is irreducible for all $s\in\F_q^\times$.
\end{lemma}

\begin{proof} See \cite[Lemma 2]{Carlos2} or \cite[Lemma 2.2]{Tung1}. \end{proof}

\medskip
We call $\lam_{\al,s}^U$ an {\it elementary} character of $U$ associated to $\al$. A {\em basic} set $D$  is a nonempty subset of $\Sig^+$ in which none of roots are on the same row or column. For each basic set $D$, we define $E(D):=\bigoplus_{\al\in D} \Irr(V_\al/[V_\al,V_\al])^\times$.

For each basic set $D$ and $\phi\in E(D)$, we define a {\em supercharacter}, also known as {\em basic} character in \cite{Carlos2},
\[\xi_{D,\phi}:=\bigotimes_{\lam_{\al,s}\in \phi}{\lam_{\al,s}}^U.\]
It turns out that each supercharacter $\xi_{D,\phi}$ is induced from a linear character of a pattern subgroup.

\begin{definition}
We define $V_D:=\bigcap_{\al\in D} V_\al$ and $\displaystyle{\lam_D:=\bigotimes_{\lam_{\al,s}\in \phi}\lam_{\al,s}|_{V_D}}$.
\end{definition}

\begin{lemma} \label{lin_lem}
We have $\xi_{D,\phi}=\lam_D^U$.
\end{lemma}

\begin{proof} See \cite[Lemma 2.5]{Tung1}. \end{proof}

\medskip
It is easy to see that $V_D$ is generated by all $X_\beta$ where $\beta\in \Sig^+-(\bigcup_{\al\in D}arm(\al))$, and $\lambda_D$ is a linear character of $V_D$.
For each basic set $D$, it can be proven that the diagonal subgroup of $\GL_n(q)$ acts transitively on $E(D)$ by conjugation. So it makes sense when we write $\lam_D$ here instead of $\lam_{D,\phi}$, and it also says that the decomposition of $\xi_{D,\phi}$ is dependent only on $D$. To know more about supercharacters, see  \cite{Is-Di,Diaconis-Thiem}. Here, we recall the main role of supercharacters as a partition of $\Irr(U)^\times$.

\begin{theorem} \label{irr_partition}
For each $\chi\in \Irr(U)^\times$, there exist uniquely a basic set $D$ and $\phi\in E(D)$ such that $\chi$ is an irreducible constituent of  $\xi_{D,\phi}$.
\end{theorem}

\begin{proof} See \cite[Theorem 1]{Carlos2} or \cite[Theorem 2.6]{Tung1}. \end{proof}

\medskip
Denote by $\Irr(\xi_{D,\phi})$ the set of all irreducible constituents of $\xi_{D,\phi}$. Here, to prove Higman's Conjecture, it suffices to prove that $|\Irr(\xi_{D,\phi})|$ is a polynomial in $q$.

Now for each basic set $D$ of size $k=|D|$, we define an associated monomial $(k\times k)$-matrix $w_D\in \GL_k(q)$. First of all, we define two partial orders on $\Sig^+$.

\begin{definition}
We define $<_r$ and $<_b$ on $\Sig^+$ as follows
\begin{itemize}
\item[(i)] $\al_{i,j}<_r \al_{l,k}$ if $j<k$ (i.e. to the right)

\item[(ii)] $\al_{i,j}<_b \al_{l,k}$ if $i<l$ (i.e. to the bottom).
\end{itemize}
\end{definition}

%\noindent
An easy way to understand these two orders is $<_r$ standing for left to right and $<_b$ for top to bottom. Notice that on a basic set, $<_r$ and $<_b$ are total orders.

Now we  fix a basic set $D$ of size $k$ ascending order of $<_r$. Let $D:=\{\tau_1,..,\tau_k\}$ where $\tau_i<_r \tau_j$ if $i<j$.
We define $w_D:=(a_{i,j})\in \GL_k(q)$ as follows
\begin{center}
$a_{i,j}:=\left\{\begin{array}{ll} 1 &\mbox{ if } \tau_j \mbox{ is the $i$-th element of $D$ in ascending order $<_b$,}\\
0&\mbox{ otherwise.}\end{array}\right.$
\end{center}
For example, if $D:=\{\al_{2,3},\al_{1,4},\al_{3,5}\}$, $|D|=3$,
\[
 \begin{array}{cccccc}
\cline{2-6}
& \multicolumn{1}{|c|}{\phantom{abc} }&\multicolumn{1}{|c|}{\phantom{abc}  }& \multicolumn{1}{|c|}{\phantom{abc}  }& \multicolumn{1}{|c|}{\alpha_{1,4}}&\multicolumn{1}{|c|}{ }  \\
\cline{2-6}
  &  & \multicolumn{1}{|c|}{ } &\multicolumn{1}{|c|}{\alpha_{2,3}} & \multicolumn{1}{|c|}{ }&\multicolumn{1}{|c|}{ }   \\
\cline{3-6}
 &   &  & \multicolumn{1}{|c|}{ } & \multicolumn{1}{|c|}{ }&\multicolumn{1}{|c|}{\alpha_{3,5}}  \\
 \cline{4-6}
  &   &   &   & \multicolumn{1}{|c|}{ }&\multicolumn{1}{|c|}{ } \\
  \cline{5-6}
  &   &   &   & &\multicolumn{1}{|c|}{ } \\
  \cline{6-6}
\end{array}
 \mbox{ ~then } w_D=\left(\begin{array}{lll}0&1&0\\1&0&0\\0&0&1\end{array}\right).
\]
%\end{small}
%
%

It is clear that $w_D$ is a monomial matrix in the Weyl group $S_k$ of $\GL_k(q)$. Here,  $w_D$ somehow gives pivots of $D$ by considering only rows and columns containing roots in $D$. Hence, it is equivalent to apply the (total) orders $<_r$, $<_b$ to these monomial matrices on their nonzero entries.

\begin{figure}
\epsfig{file=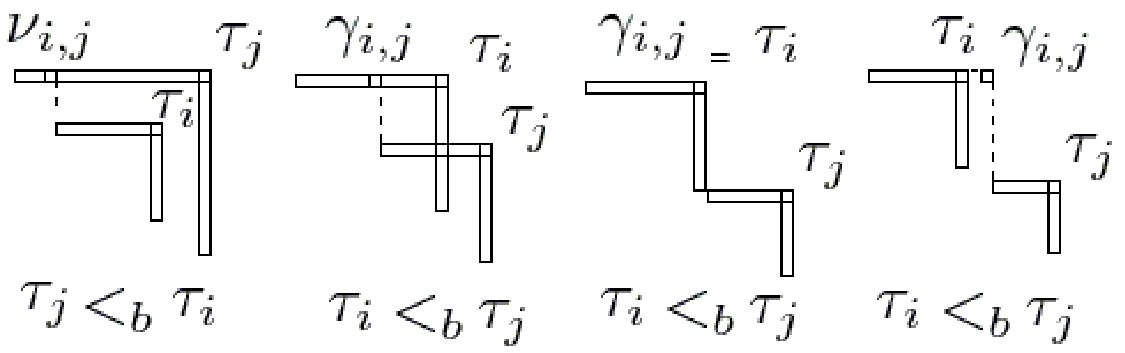,width=3.5in}
\caption{\sl Positions of $\nu_{i,j}$ and $\gamma_{i,j}$}
\label{fig:PosPic}
\end{figure}

For each pair $0<i<j\leq k$, if $\tau_i<_b\tau_j$, let $\gamma_{i,j}$ be the root on the row of $\tau_i$ such that $\gamma_{i,j}+\tau_j\in \Sigma^+$; otherwise, i.e. $\tau_j<_b\tau_i$, let $\nu_{i,j}$ be the root on the row of $\tau_j$ such that $\nu_{i,j}+\tau_i\in\Sig^+$. For example, $\tau_i:=\al_{m,i}$, $\tau_j:=\al_{l,j}$ where $i<j$, so if $\al_{m,i}<_b\al_{l,j}$, i.e. $m<l$, then $\gamma_{i,j}=\al_{m,l-1}$; otherwise, if $\al_{l,j}<_b\al_{m,i}$, i.e. $l<m$, then $\nu_{i,j}=\al_{l,m-1}$. It is easy to see that $\nu_{i,j}$ exists if and only if two hooks $h(\tau_i)$ and $h(\tau_j)$ are parallel, otherwise $\gamma_{i,j}$ exists (Figure \ref{fig:PosPic}).

Let $\Gamma_D$ be the set of all $\gamma_{i,j}$, let $\Lambda_D$ be the set of all $\nu_{i,j}$, and $\Delta_D:=\Gamma_D\cup\Lambda_D$. By definitions for the existence of $\gamma_{i,j}$ and $\nu_{i,j}$, we have $\Gamma_D\cap \Lambda_D=\emptyset$.

\begin{definition}
We define $R_D:=\la X_\gamma : \gamma\in \Gamma_D\ra$ and $C_D:=\la X_\nu : \nu\in \Lambda_D\ra$.
\end{definition}

The next lemma provides interesting correspondences between the size of $D$ and $\Delta_D$, and between $w_D$ and $\Gamma_D$, $\Lambda_D$. Moreover, it shows that $\la V_D,R_D\ra=V_D R_D$, and the pattern subgroups $R_D$, $C_D$ are only determined by $w_D$ in a natural way.

\begin{lemma} \label{main_lemma}
Let $D$ be a basic set of size $k$. The following are true.
\begin{itemize}
\item[(i)] $\Delta_D$ is closed and %%$\la R_D,C_D\ra=
$\la X_\al : \al\in \Delta_D\ra$ is isomorphic to $U_k(q)$.

\item[(ii)] $\Gamma_D$ is closed. For each pair $i<j$, if $\gamma_{i,s},\gamma_{j,r}$ exist and $\gamma_{i,s}+\gamma_{j,r}\in \Sig^+$, then $s=j$ and $\gamma_{i,j}+\gamma_{j,r}=\gamma_{i,r}$.

\item[(iii)] $\Lambda_D$ is closed. For each pair $i<j$, if $\nu_{i,s},\nu_{j,r}$ exist and $\nu_{i,s}+\nu_{j,r}\in \Sig^+$, then $s=j$ and $\nu_{i,j}+\nu_{j,r}=\nu_{i,r}$.

\item[(iv)] $R_D$ is isomorphic to $U_k(q)\cap {}^{w_D}U_k(q)$ and \\$C_D$ is isomorphic to $U_k(q)\cap {}^{w_0w_D}U_k(q)$ where $w_0$ is the longest element in the Weyl group $S_k$ of $\GL_k(q)$, also known as
    \[w_0=\begin{small}\left(\begin{array}{ccc}0&\cdots &1\\ \vdots& \iddots & \vdots\\ 1& \cdots & 0 \end{array}\right)\end{small}.\]

\item[(v)] $V_D R_D$ is a pattern subgroup of $U$ and $R_D$ normalizes $V_D$.

\end{itemize}

\end{lemma}

For example, let $D:=\{\al_{1,2},\al_{3,4},\al_{4,5},\al_{2,6}\}$ be a basic set in $\Sig_6^+$.
\begin{center}
\begin{small}
\[U_7(q)= \begin{array}{cccccc}
\cline{1-6}       %line 1
 \multicolumn{1}{|c|}{\phantom{\al_1}} & \multicolumn{1}{|c|}{  \al_{1,2} } & \multicolumn{1}{|c|}{ } & \multicolumn{1}{|c|}{  }&\multicolumn{1}{|c|}{ } & \multicolumn{1}{|c|}{ } \\

\cline{1-6}    %line 2
 & \multicolumn{1}{|c|}{\phantom{\al_2}} & \multicolumn{1}{|c| }{ } & \multicolumn{1}{|c|}{ } & \multicolumn{1}{|c|}{   } & \multicolumn{1}{|c|}{\al_{2,6} }\\

\cline{2-6}   %line 3
 &  & \multicolumn{1}{|c|}{\phantom{\al_3}} & \multicolumn{1}{|c|}{\al_{3,4} } & \multicolumn{1}{|c|}{ } & \multicolumn{1}{|c|}{ }  \\

\cline{3-6}  %line 4
 &   &   & \multicolumn{1}{|c|}{\phantom{\al_4} } & \multicolumn{1}{|c|}{\al_{4,5} } & \multicolumn{1}{|c|}{ }  \\

\cline{4-6}             %line 5
 &   &   & &\multicolumn{1}{|c|}{\phantom{\al_5}} & \multicolumn{1}{|c|}{   }  \\

\cline{5-6}              %line 6
 &   &   & & &\multicolumn{1}{|c|}{\phantom{\al_6}} \\

\cline{6-6}               %line 7
\end{array}
\]
\end{small}
and
\begin{small}
\[R_D= \begin{array}{ccc}
\cline{1-3}       %line 1
 \multicolumn{1}{|c|}{\phantom{\al_1}} & \multicolumn{1}{|c|}{  \al_{1,2} } & \multicolumn{1}{|c|}{ } \\

\cline{1-3}    %line 2
 & \multicolumn{1}{ c }{\phantom{\al_2}} & \multicolumn{1}{ c  }{ } \\

\cline{3-3}   %line 3
 &  & \multicolumn{1}{ |c| }{\phantom{\al_3}}  \\

\cline{3-3}  %line 4

\end{array}\ \ , \ \ \ \
C_D= \begin{array}{ccc}
%\cline{1-3}       %line 1
 \multicolumn{1}{ c }{\phantom{\al_1}} & \multicolumn{1}{ c }{ } & \multicolumn{1}{  c }{ } \\

\cline{2-3}    %line 2
 & \multicolumn{1}{ |c| }{\phantom{\al_2}} & \multicolumn{1}{ |c|  }{ } \\

\cline{2-3}   %line 3
 &  & \multicolumn{1}{  c  }{\phantom{\al_3}}  \\

%\cline{3-3}  %line 4

\end{array}\ \ .
\]
\end{small}
\end{center}

The next result is the main theorem, which provides a correspondence between supercharacters $\xi_{D,\phi}$ and pattern subgroups $R_D$.

\begin{theorem} \label{main_thm}
Let $\xi_{D,\phi}$ be a supercharacter. The following are true.
\begin{itemize}
\item[(i)]  $\xi_{D,\phi}=(\lam_D^{V_D R_D })^U$.

\item[(ii)] For each $\chi\in \Irr(V_D R_D,\lam_D)$, $\chi^U\in \Irr(\xi_{D,\phi})$.

\item[(iii)]If $\chi_1\neq\chi_2\in \Irr( V_D R_D ,\lam_D)$, then ${\chi_1}^U\neq {\chi_2}^U$.

\end{itemize}
\end{theorem}

Therefore, to decompose $\xi_{D,\phi}$, it suffices to decompose $\lam_D^{V_D R_D}$. Moreover, the induced character $\lambda_D^{V_D R_D}$ is equal to $({\lambda_D|_{V_D\cap R_D}}^{R_D})_{V_D R_D}\otimes \theta$, where $\theta$ is some linear character of $V_D R_D$ (in Lemma \ref{sub_lemma}). We see that ${\lambda_D|_{V_D\cap R_D}}^{R_D}$ is a `very special' constituent of the regular character $1^{R_D}$. Hence, the decomposition method of all supercharacters $\xi_{D,\phi}$ of $U_n(q)$ with the same $w_D$ is generally restricted to the one of the regular character $1^{R_D}$.

Here, we attempt to make a link for this special pattern $R_D=U_k(q)\cap {}^{w_D}U_k(q)$ in Lemma \ref{main_lemma}. Denote $U\cap {}^w U$ by $U_w$, where $U=U_n(q)$ and $w\in S_n$ is the Weyl group of $\GL_n(q)$, Thompson \cite{Thompson2} conjectured that for each pair $r,s\in S_n$, the cardinality of the double coset $U_r \backslash U / U_s$ is a polynomial in $q$ with integer coefficients. In addition, $U_w$ also takes an important role when one studies $\GL_n(q)$ as groups with $(B,N)$-pair, such as the Bruhat decomposition.

From Theorem \ref{main_thm} and Lemma \ref{main_lemma} (v), we obtain a nice decomposition of $\xi_{D,\phi}$.

\begin{corollary} \label{cor:super-decomp}
Let $\xi_{D,\phi}$ be a supercharacter. The following are tru:.
\begin{itemize}
\item[(i)] $\Irr(\xi_{D,\phi})=\{\chi^U:\, \chi\in \Irr(V_DR_D,\lam_D)\}$,

\item[(ii)] $\xi_{D,\phi}=\sum_{\chi\in \Irr(V_DR_D,\lam_D)}\chi(1)\chi^U$.
\end{itemize}

\end{corollary}

Theorem \ref{irr_partition}, Lemma \ref{main_lemma} and Corollary \ref{cor:super-decomp} give a clear proof for the following corollary, which is a different version of \cite[Theorem 1.4]{Carlos2}.

\begin{corollary}
$(\xi_{D,\phi},\xi'_{D',\phi'})=\left\{\begin{array}{ll} |V_DR_D:V_D|&\mbox{ if } (D,\phi)=(D',\phi'),\\0&\mbox{ otherwise.}\end{array}\right.$
\end{corollary}

As an application to $U_{13}(q)$, we answer the conjecture by Isaacs-Karagueuzian \cite{Is-Kar} that $U_{13}(2)$ has a unique pair of irrational irreducible characters of degree $2^{16}$. This conjecture is solved with an affirmative answer by  Marberg \cite{Marberg} and generalized by Evseev \cite{Anton1}. Here, we work independently to obtain representations and constructions of all irreducible constituents of the supercharacter of $U_{13}(q)$, which gives the irrational pair when $q=2$.  With the definition that a character is {\em well-induced} if it is induced from a linear character of some pattern subgroup by Evseev \cite{Anton1}, these two irrational characters are not well-induced. Thus, it provides a more explainable script to the generalization of not well-induced characters.  Finally, we list two families of supercharacters $\xi_{D,\phi}$ which have exactly one irreducible constituent, i.e. $\xi_{D,\phi}=m\ \chi$ for some $\chi\in \Irr(U)$.

%%%%%%%%%%%%%%%%%%%%%%%%%%%%%%%%%%%%%%%%%%%%%%%%%%%%%%%%%%%%%%%%%%%%%%%%%%%%%%%%%%%%%%%%%%

\section{All proofs}
\label{sec:proofs}

Here we mainly prove Theorem \ref{main_thm} to give a natural correspondence between supercharacters $\xi_{D,\phi}$ and pattern subgroups $U_k(q)\cap {}^{w_D}U_k(q)$ where $k=|D|$. First, we shall prove Lemma \ref{main_lemma}.

%%%%%%%%%%%%%%%%%%%%%%%%%% PROOF OF MAIN LEMMA *********************************

%\bigskip
\begin{proof}[Proof of Lemma \ref{main_lemma}] \label{proof:lem2.7}
Suppose that $D:=\{\tau_1,...,\tau_k\}$ in ascending order $<_r$.

(i) If we rearrange $D$ in ascending order of $<_b$ to be $\{\theta_1,..,\theta_k\}$, then on the row of $\theta_i$, $\Delta_D$ has $(k-i)$ roots and the row of $\theta_k$ does not have any root in $\Delta_D$.

For each pair $i<j\in[1,k]$, let $\omega_{i,j}\in \Delta_D$ be the root on the row of $\tau_i$ such that $\omega_{i,j}+\tau_j\in \Sig^+$. (Notice that $\omega_{i,j}$ is either $\gamma\in \Gamma_D$ or $\nu\in\Lambda_D$.) Hence, if $\tau_i=\al_{i_1,i_2}<_b \tau_j=\al_{j_1,j_2}$, i.e. $i_1<j_1$, we have $\omega_{i,j}=\al_{i_1,j_1-1}$. Therefore, for each $\omega_{i,j}=\al_{i_1,j_1-1}<_r\omega_{m,l}=\al_{m_1,l_1-1}\in\Delta_D$, if $\omega_{i,j}+\omega_{m,l}\in\Sig^+$, then it must be $j_1=m_1$, and $\omega_{i,j}+\omega_{j,l}=\al_{i_1,l_1-1}=\omega_{i,l}$. This shows that $\Delta_D$ is closed, and the longest root in $\Delta_D$ is $\omega_{1,2}+...+\omega_{k-1,k}=\omega_{1,k}$. So $\omega_{i,j}$ corresponds to $\al_{i,j-1}$ in the positive root set $\Sig_{k-1}^+$. Therefore, $\la X_\al:\al\in \Delta_D\ra$ is a pattern subgroup isomorphic to $U_k(q)$.

\medskip
(ii) With the same argument as in (i), by definitions of $\gamma_{i,s}$, $\gamma_{j,r}$, if $\gamma_{i,s}+\gamma_{j,r}\in \Sig^+$, then $s=j$.
By the transitivity of $<_r,<_b$ on $\tau_i,\tau_j,\tau_r$ , from $\tau_i<_r,<_b \tau_j$ and $\tau_j <_r,<_b \tau_{k}$, we have $\tau_i <_r,<_b \tau_r$. So $\gamma_{i,r}$ exists and  $\gamma_{i,j}+\gamma_{j,r}=\gamma_{i,r}$ follows.

\medskip
(iii)  The same argument of (ii) holds for $\nu_{i,s}$ and $\nu_{j,r}\in\Lambda_D$.

\medskip
(iv) Let $w_D{:=}(w_{i,j})\in S_k \subset \GL_k(q)$.  Since $w_D$ is a monomial matrix,  $w_D^{-1}=w_D^T$, the transpose of $w_D$. For each $X=(x_{i,j})\in U_k(q)$, we observe $Y:=w_D\cdot X\cdot {w_D}^{-1}$. Let $Y=(y_{i,j})$. For each pair $i<j$, we have
  $y_{i,j}=\sum_{s,r\in[1,k]}w_{i,s}x_{s,r}w_{j,r}$.
Since $i,j$ are fixed, there exist uniquely $1\leq f,h\leq k$ such that $w_{i,f}=1=w_{j,h}$, others $w_{i,s}=0=w_{j,r}$. Hence, $y_{i,j}=w_{i,f}x_{f,h}w_{j,h}$.

Since $h\neq f$ and all $x_{s,r}=0$ if $r<s$, we have the following
  \begin{itemize}
  \item $y_{i,j}=0$  if $f>h$, i.e. $w_{i,f}<_b w_{j,h}$ and $w_{j,h}<_r w_{i,f}$;
  \item $y_{i,j}\neq 0$ if $f<h$, i.e. $w_{i,f}<_b w_{j,h}$ and $w_{i,f}<_r w_{j,h}$.
  \end{itemize}
So $R_D$ is isomorphic to $U_k(q)\cap {}^{w_D}U_k(q)$ by the definition of $\gamma_{i,j}\in\Gamma_D$. And, hence, $C_D$ is isomorphic to $U_k(q)\cap {}^{w_0\cdot w_D}U_k(q)$ by (i), (ii), (iii) and $\Delta_D=\Gamma_D\cup\Lambda_D$.

\medskip
(v) From the definition of $\gamma_{i,j}$, it is easy to check that $R_D$ normalizes $V_D$. Hence, $V_D R_D$ is a pattern subgroup of $U$.
\end{proof}

%%%%%%%%%%%%%%%%%%%%%%%%%%%%% DONE FOR MAIN LEMMA %%%%%%%%%%%%%%%%%%%%%%%%%%%%%%%%%

\medskip
Set $K_D:=\la X_\al: X_\al\subset V_D \mbox{ and } \al\notin D\ra=\la X_\al: X_\al\subset V_D\cap \ker(\lambda_D)\ra$. It is clear that $K_D$ is normal in $V_D$, $[V_D:K_D]=q^{|D|}$, and $V_D=K_D\cdot \prod_{\tau\in D}X_\tau$. To prove Theorem \ref{main_thm}, we need the following lemma.
% For each character $\chi$ of a group $G$, denote $Z(\chi):=\{x\in G:~ |\chi(x)|=\chi(1)\}$.

%%%%%%%%%%%%%%%%%%%%%%       SUB_LEMMA  %%%%%%%%%%%%%%%%%%%%%%%%%%%

\begin{lemma} \label{sub_lemma}

Let $\xi_{D,\phi}$ be a supercharacter. The following are true.
\begin{itemize}

\item[(i)] $K_D{\subset} \ker(\lam_D^{V_DR_D})$. Moreover, $\lam_D^{V_DR_D}(x){=}|V_DR_D{:}V_D|\lam_D(x)$ for all $x\in V_D$.

\item[(ii)] $(K_D\cap R_D)\unlhd R_D$ and $(V_D\cap R_D)/(K_D\cap R_D)\subset Z(R_D/(K_D\cap R_D))$.

\item[(iii)] Let $\overline{\phi}_D:=\{\lam_{\al,s}\in \phi:\, X_\al\nsubseteq R_D\}$. We have
\[\lam_D^{V_D R_D}=({\lam_D|_{V_D\cap R_D}}^{R_D})_{V_D R_D}\otimes \left(\bigotimes_{\lam_{\al,s}\in \overline{\phi}_D}(\lam_{\al,s}|_{V_D})_{V_DR_D}\right).\]

\end{itemize}

\end{lemma}

\begin{proof} (i)   It suffices to show the statement for all $X_\al\subset V_D$. By Lemma \ref{main_lemma} (v) $V_D\unlhd V_D R_D$, for all $x\in V_D$
we have
\[\lam_D^{V_DR_D}(x)=\frac{1}{|V_D|}\sum_{y\in V_DR_D}\lam_D(x^y).\]
For each $x\in X_\al$, we suppose that there is $X_\beta\subset V_DR_D$ such that $\al+\beta\in \Sig^+$, hence $X_{\al+\beta}\subset V_D$. We shall show that $\lam_D(x^y)=\lam_D(x)$ for all $y\in X_\beta$.

Since $X_\tau\cap [V_D,V_D]=\{1\}$ for all $\tau\in D$, we have $X_{\al+\beta}\subset K_D\subset \ker(\lam_D)$.
Thus, $[\lam_D(x),\lam_D(y)]=\lam_D([x,y])=1$ since $[x,y]\in X_{\al+\beta}$, i.e. $\lam_D(x)^{-1}\lam_D(x^y)=1$.

\medskip
(ii) By the definition of $K_D\unlhd V_D$ and $V_D=K_D\cdot \prod_{\tau\in D}X_\tau$, it suffices to show that $(K_D\cap R_D)\unlhd R_D$. This is clear because for all $X_\al\subset K_D\cap R_D$ and all $X_\beta\subset R_D$, either $\al+\beta\not\in \Sig^+$ or $X_{\al+\beta}\subset K_D\cap R_D$.

\medskip
(iii) The inflations to $V_D R_D$ of ${\lam_D|_{V_D\cap R_D}}^{R_D}$ and $\lam_{\al,s}|_{V_D}$, for all $\lam_{\al,s}\in \overline{\phi}_D$,  follow directly from (i).
\end{proof}

%%%%%%%%%%%%%% DONE FOR SUB_LEMMA %%%%%%%%%%%%%%%%%%%%%%%%%%%%%%%%%%%

\begin{remark}
By Lemma \ref{sub_lemma} (iii), if $R_D\cap V_D=\{1\}$, then $\lam_D^{V_DR_D}$ is equivalent to the regular character $1^{R_D}$ of $R_D$. In general, $\lam_D^{V_DR_D}$ is equivalent to a constituent of $1^{R_D}$ with $R_D\cap K_D$ in the kernel. Now we prove Theorem \ref{main_thm}.
\end{remark}

%%%%%%%%%%%%%%%%PROOF OF MAIN THEOREM **********************

\smallskip
\begin{proof}[Proof of Theorem \ref{main_thm}]
(i) This is clear by the transitivity of induction.

\medskip
(ii) We suppose that $D:=\{\tau_1,...,\tau_k\}$ in ascending order $<_r$ and for $s_i\in\F_q^\times$ let $\lam_D:=\bigotimes_{\tau_i\in D}\lam_{\tau_i,s_i}|_{V_D}$.

First, we show that for each $\chi\in \Irr(V_DR_D,\lam_D)$, $\chi^U$ is irreducible. By the transitivity of induction, we shall induce $\chi$ from $V_DR_D$ to $U$ by a sequence of inductions along the arms of $\tau_1,\tau_2,...,\tau_k$ respectively by $<_r$ order. Now we setup these such induction steps.

For each $\tau_i\in D$, let $A(\tau_i)=\{\al \in arm(\tau_i): X_\al\nsubseteq V_DR_D\}$, and $c_i=|A(\tau_i)|$. Let $d_0=0$ and $d_i=d_{i-1}+c_i$ for all $i\in[1,k]$. Now if $c_i>0$, $i\in[1,k]$, we arrange $A(\tau_i)$ in decreasing order $<_r$ to be $\{\beta_{d_{i-1}+1},...,\beta_{d_{i-1}+c_i}\}$. Let $M_0=V_DR_D$, $M_{i+1}=M_i\rtimes X_{\beta_i}$ for all $i\in[1,d_k]$. It is clear that $M_{d_k+1}=U$ and $X_{\beta_j}$ normalizes $M_j$; hence, this sequence of pattern subgroups is well-defined.

For each $\beta_j\in arm(\tau_i)$, $j\in[1,d_k]$, there exists uniquely $\delta\in leg(\tau_i)$ such that $\beta_j+\delta=\tau_i$ and $X_\delta\subset K_D$, since if $X_\delta\nsubseteq K_D$, there exists $\tau_m\in D$ such that $\delta\in arm(\tau_m)$, so $\tau_i<_r\tau_m$, $\tau_i<_b\tau_m$, and this implies $\beta_j=\gamma_{i,m}$. We number this $\delta$ as $\delta_j$, and let $L(D):=\{\delta_j:~j\in [1,d_k]\}$. By Lemma \ref{sub_lemma} (i), $X_\delta\subset \ker(\chi)$ for all $\delta\in L(D)$.
Now we proceed the induction of $\chi$ from $V_DR_D$ to $U$ via a sequence of pattern subgroups along the arms of all $\tau_i\in D$, namely from $M_0$ to $M_1$, ..., $M_{d_k+1}=U$.

Suppose that $\chi^{M_j}\in \Irr(M_j)$ for some $M_j$, $j\in[1,d_k+1]$ and $X_{\delta_t}\subset \ker(\chi^L)$ for all $t\in[j,d_k]$. If $j=d_k+1$, the process is complete. Otherwise, the next induction step is from $M_j$ to $M_{j+1}=M_j X_{\beta_j}$, and we suppose that it happens on the arm of $\tau_i$. For each $x\in X_{\beta_j}^\times$, since $[X_{\delta_j},x]= X_{\tau_i}$, there is some $y\in X_{\delta_j}$ such that $\lam_{\tau_i,s_i}([y,x])\neq 1$ and
\begin{center}
${}^x(\chi^{M_j})(y)=\chi^{M_j}(y^x)=\chi^{M_j}([y,x]y)=\lam_{\tau_i,s_i}([y,x])\chi^{M_j}(y)\neq \chi^{M_j}(y)=\chi^{M_j}(1)$.
\end{center}
Hence, $X_{\delta_j}\nsubseteq \ker({}^x(\chi^{M_j}))$, and ${}^x(\chi^{M_j})\neq \chi^{M_j}$ for all $x\in X_{\beta_j}^\times$. This shows that the inertia group $I_{M_jX_{\beta_j}}(\chi)=M_j$ and $\chi^{M_jX_{\beta_j}}\in \Irr(M_jX_{\beta_j},\lambda_D)$.

It is easy to check directly that $X_{\delta_t}\subset \ker(\chi^{M_jX_{\beta_j}})$ for all $t\in [j+1,d_k]$ by using $[X_{\beta_j},X_{\delta_t}]\subset \ker(\chi^{M_j})$. Therefore, $\chi^U$ is irreducible for all $\chi\in \Irr(V_DR_D,\lambda_D)$ by induction on $j$.

\medskip
(iii) Now suppose that $\chi_1\neq \chi_2\in \Irr(V_DR_D,\lambda_D)$ and ${\chi_1}^{M_j}\neq {\chi_2}^{M_j}$ for some $M_j$. As above, it is enough to show that $\chi_1^{M_jX_{\beta_j}}\neq \chi_2^{M_jX_{\beta_j}}$, where $\beta_j\in arm(\tau_i)$.
Notice that $X_{\delta_j}\subset \ker({\chi_1}^{M_j})\cap \ker({\chi_2}^{M_j})$.

By Mackey formula with the double coset $M_j\backslash M_jX_{\beta_j}/M_j$ represented by $X_{\beta_j}$, $({\chi_1}^{M_jX_{\beta_j}},{\chi_2}^{M_jX_{\beta_j}})=\sum_{x\in X_{\beta_j}}({\chi_1}^{M_j},{}^x({\chi_2}^{M_j}))$. Using the same argument as in (ii), we have $X_{\delta_j}\nsubseteq \ker({}^x({\chi_2}^{M_j}))$ for all $x\in X_{\beta_j}^\times$.
Hence, ${}^x({\chi_2}^{M_j})\neq {\chi_1}^{M_j}$ for all $x\in X_{\beta_j}^\times$ since $X_{\delta_j}\subset \ker({\chi_1}^{M_j})$. Thus, $({\chi_1}^{M_jX_{\beta_j}},{\chi_2}^{M_jX_{\beta_j}})=({\chi_1}^{M_j},{\chi_2}^{M_j})=0$ since ${\chi_1}^{M_j}\neq {\chi_2}^{M_j}$ by the above assumption on $M_j$.
\end{proof}

%%%%%%%%%%%%%%%%%%%%%%Done the proof of Main Theorem %%%%%%%%%%%%%%%%%%%%%

\medskip
Notice that  $V_DR_D$ is not normal in $U$. In the proof of Theorem \ref{main_thm}, although all inductions from $V_DR_D$ to $U$ are irreducible, Clifford correspondence can not be applied. The technique of a sequence of inductions from $M_j$ to $M_{j+1}\subset N_U(M_j)$ has been used to control distinct induced characters.

Since $V_D$ is normal in $V_DR_D$ and $V_DR_D/V_D\cong R_D/(V_D\cap R_D)$, by Theorem \ref{main_thm} and Lemma \ref{sub_lemma} (iii), we only need to decompose  ${\lam_D|_{V_D\cap R_D}}^{R_D}$ instead of decomposing the supercharacter $\xi_{D,\phi}=\lam_D^U$. Hence, all work is restricted to a pattern subgroup of $U_k(q)$ where $k=|D|<n$.

%%%%%%%%%%%%%%%%%%%%%%%%%% PROOF OF SUPER_DECOMPOSITION %%%%%%%%%%%%%%%%%%%%%%%%%%%%%%%%

\medskip
\begin{proof}[Proof of Corollary \ref{cor:super-decomp}]
Theorem \ref{main_thm} provides a one-to-one correspondence on the multiplicities and degrees between two constituent sets $\Irr(V_DR_D,\lam_D)$ and $\Irr(\xi_{D,\phi})$,  i.e. $|\Irr(V_DR_D,\lam_D)|=|\Irr(\xi_{D,\phi})|$, and if $\chi\in \Irr(V_DR_D,\lam_D)$ has multiplicity $t$ then $\chi^U\in \Irr(\xi_{D,\phi})$ also has multiplicity  $t$, and $\chi^U(1)=|U:V_DR_D|\chi(1)$.
Therefore, it is enough to show that $\chi\in \Irr(R_D,\lam_D|_{V_D\cap R_D})$ has multiplicity $\chi(1)$.

By Lemma \ref{sub_lemma} (i), $K_D\cap V_D\subset \ker(\lam_D|_{V_D\cap R_D}) \cap \ker({\lam_D|_{V_D\cap R_D}}^{R_D})$ is normal in $R_D$. So $\lam_D|_{V_D\cap R_D}$ can be considered as a linear character of the quotient group $R_D/(K_D\cap R_D)$. By Lemma \ref{sub_lemma} (ii), $(V_D\cap R_D)/(K_D\cap R_D)\subset Z(R_D/(K_D\cap R_D))$,  $\lam_D|_{V_D\cap R_D}$ is a linear character of the center, and the claim holds.
\end{proof}

% Say a little about action of split torus (transitive) to obtain permutation character $1^{R_D}$ if $V_D\cap R_D=\{1\}$, otherwise it is a difference of permutation characters., i.e. exist $X_\tau\subset R_D$, then $\xi_{D-\tau}=\lam_{V_{D-\tau}}+\sum_{t\in T_\tau\cong \F_q^\times} {}^t\lam_D$ where $\xi_{D-\tau}\in \Irr(K_D\prod_{\al\in D-\tau}X_\al)$ such that $\xi_{D-\tau}|_{X_\al}=\lam_D|_{X_\al}$ for all $\al\in D-\tau$

%%%%%%%%%%%%%%%%%%%%%%%%%%% APLLICATIONS****************************

\section{Applications}
\label{sec:app}

%\noindent {\bf Example 1:}
\begin{example}
Here we list two families of supercharacters $\xi_{D,\phi}$ of $U_n(q)$ which have only one irreducible constituent, i.e. $|\Irr(\xi_{D,\phi})|=1$. Without loss of generality, we suppose that $\al_{1,-}\in D$.

$D_1:=\{\al_{1,k},\al_{2,2k-1},\al_{3,2k-2},...,\al_{k,k+1},\al_{k+1,2k}\}$ where $4\leq2k<n$. We have
\[w_{D_1}=\left(\begin{small} \begin{array}{ccc} 1&.&.\\.&w&.\\.&.&1\end{array}\end{small}\right)\]
where $w$ is the longest element in the Weyl group $S_{k-1}$ of $\GL_{k-1}(q)$.

$D_2:=\{\al_{1,2},\al_{2,3},...,\al_{2m-1,2m}\}$ where $2\leq2m<n$, which gives $w_{D_2}$ equal to the identity $(2m-1){\times} (2m-1)$-matrix $I_{2m-1}$.

By Lemma \ref{main_lemma}, $R_{D_1}$ is the largest special subgroup $q^{1+2(k-1)}$ in $U_{k+1}(q)$, i.e. $[R_{D_1},R_{D_1}]=Z(R_{D_1})=\Phi(R_{D_1})$ the Frattini subgroup of $R_{D_1}$.  It is known that a special subgroup of type $q^{1+2t}$ has $q^{2t}$ linear characters and $q-1$ almost faithful irreducible characters of degree $q^{t}$, see \cite[Corollary 2.3]{Tung1}. (Recall that $\chi\in \Irr(G)$ is {\em almost faithful} if $Z(G)\not\subset \ker(\chi)$.) Since $V_{D_1}\cap R_{D_1}=X_{\al_{1,k}}=Z(R_{D_1})$, $\lam_{D_1}|_{X_{\al_{1,k}}}\neq 1_{X_{\al_{1,k}}}$. Hence, $\lam_{D_1}^{V_{D_1}R_{D_1}}$ has only one irreducible constituent of degree $q^{k-1}$ with multiplicity $q^{k-1}$. By Corollary \ref{cor:super-decomp}, $\xi_{D_1,\phi}$ has only one constituent of degree $q^{(k-1)+|U:V_{D_1}R_{D_1}|}$ with multiplicity $q^{k-1}$.

We decompose $\xi_{D_2,\phi}$ by induction on $m$.  By Corollary \ref{cor:super-decomp} and Lemma \ref{sub_lemma}, we have $|\Irr(\xi_{D_2,\phi})|=|\Irr(R_{D_2},\lam_{D_2}|_{V_{D_2}\cap R_{D_2}})|$. Let $D'_2:=\{\al_{1,2},\al_{2,3}$, ..., $\al_{2m-3,2m-2}\}$. Since $w_{D_2}=I_{2m-1}$, by Lemma \ref{main_lemma}, $R_{D_2}\simeq U_{2m-1}(q)$.
 It is easy to check that $D'_2=\{\al\in D: X_\al\subset R_{D_2}\}$ and ${\lam_{D_2}|_{V_{D_2}\cap R_{D_2}}}^{R_{D_2}}=q\,\xi_{D'_2,\phi'}$ where $\xi_{D'_2,\phi'}$ is a supercharacter of $U_{2m-1}(q)$ with $\phi'=\phi-\{\lambda_{\al_{2m-2,2m-1},s_{2m-2}},\lambda_{\al_{2m-1,2m},s_{2m-1}}\}$.  Hence, by the hypothesis of induction on $m$, it suffices to check $D_2=\{\al_{1,2},\al_{2,3},\al_{3,4}\}$, which is $D_1$ with $k=2$.
\end{example}

%\medskip
%\noindent {\bf Example 2:}
\begin{example}
Isaacs and Karagueuzian \cite{Is-Kar} conjectured that $U_{13}(2)$ has a unique pair of irrational irreducible characters of degree $2^{16}$. This conjecture was answered by Marberg (\cite[Theorem 9.2]{Marberg}). By constructing decomposition trees, he also found out the exact supercharacter whose irreducible constituent set contains this pair. Evseev (\cite[Theorem 1.4]{Anton1}) generalized the result by applying a reduction algorithm process to obtain $q(q-1)^{13}$ irreducibles of $U_{13}(q)$ which cannot be constructed by inducing from a linear character of  some pattern subgroup. Those characters are called {\em not well-induced}. We work independently and point out the supercharacter with these properties as well.
\end{example}

To deal with extensions of linear characters from a subgroup, we use the following property whose proof is quite obvious.

\begin{lemma} \label{lem:ExtChar}
For $H\leq G$ and $\lam\in\Irr(H)$, if $[G,G]\leq \ker(\lam)$ then $\lam$ extends to $G$.
%\begin{itemize}
%\item[(i)] If $\lam$ is linear and $[G,G]\subset \ker(\lam)$ then $\lam$ extends to $G$.
%\item[(ii)] If $G\simeq K\rtimes H$ where $K\trianglelefteq G$ then $\lam$ inflates to $G$.
%\end{itemize}
\end{lemma}

Let $U{:=}U_{13}(q)$ and $D{:=}\{\al_{1,4},\al_{2,5},\al_{3,9},\al_{4,10},\al_{5,6},\al_{6,7},\al_{7,8},\al_{8,11},\al_{9,12}\}\subset \Sig_{12}^+$.  By Lemma \ref{main_lemma}, we have
\[w_D=\begin{tiny}\left(
\begin{array}{ccccccccc}
1&.&.&.&.&.&.&.&.\\
.&1&.&.&.&.&.&.&.\\
.&.&.&.&.&1&.&.&.\\
.&.&.&.&.&.&1&.&.\\
.&.&1&.&.&.&.&.&.\\
.&.&.&1&.&.&.&.&.\\
.&.&.&.&1&.&.&.&.\\
.&.&.&.&.&.&.&1&.\\
.&.&.&.&.&.&.&.&1\\
\end{array}
\right) \end{tiny}\in S_9\leq \GL_9(q)
\]
 and
\begin{small}
\[R_D= \begin{array}{cccccccc}
\cline{1-8}       %line 1
 \multicolumn{1}{|c|}{\al_1} & \multicolumn{1}{|c|}{  } & \multicolumn{1}{|c|}{  } & \multicolumn{1}{|c|}{ ~*~}&\multicolumn{1}{|c|}{\bullet} & \multicolumn{1}{|c|}{\bullet}& \multicolumn{1}{|c|}{\bullet} & \multicolumn{1}{|c|}{\bullet} \\

\cline{1-8}    %line 2
 & \multicolumn{1}{|c|}{\al_2} & \multicolumn{1}{|c| }{ } & \multicolumn{1}{|c|}{ } & \multicolumn{1}{|c|}{ * } & \multicolumn{1}{|c|}{\bullet} & \multicolumn{1}{|c|}{\bullet} &\multicolumn{1}{|c|}{\bullet}  \\

\cline{2-8}   %line 3
 &  & \multicolumn{1}{|c|}{\al_3} & \multicolumn{1}{ c }{ } & \multicolumn{1}{ c }{ } & \multicolumn{1}{ c }{ } & \multicolumn{1}{|c|}{ } & \multicolumn{1}{|c|}{ } \\

\cline{3-3}\cline{7-8}  %line 4
 &   &   & \multicolumn{1}{c }{ } & \multicolumn{1}{ c }{ } & \multicolumn{1}{ c }{ } & \multicolumn{1}{|c|}{}&\multicolumn{1}{|c|}{ } \\

\cline{5-8}             %line 5
 &   &   & &\multicolumn{1}{|c|}{\al_5} & \multicolumn{1}{|c|}{ * } & \multicolumn{1}{|c|}{\bullet} & \multicolumn{1}{|c|}{\bullet} \\

\cline{5-8}              %line 6
 &   &   & & &\multicolumn{1}{|c|}{\al_6}&\multicolumn{1}{|c|}{ * }&\multicolumn{1}{|c|}{\bullet} \\

\cline{6-8}               %line 7
 &   &   & & & &\multicolumn{1}{|c|}{\al_7}&\multicolumn{1}{|c|}{ * } \\
\cline{7-8}               %line 8
&   &   & & & & &\multicolumn{1}{|c|}{\al_8} \\
\cline{8-8}               %line 9
\end{array}\ \ {\normalsize  \leq\ U_9(q)} \ .
\]
\end{small}

As above $|U{:}V_DR_D|=q^{12}$, $\lambda_D^U(1)=q^{27}$, $R_D=U_9(q)\cap {}^{w_D}U_9(q)$,
$D\cap\Gamma_D=$ all $*$'s, $K_D\cap R_D=$ all $X_\bullet\mbox{'s}\subset \ker(\lambda_D)$, and $V_D\cap R_D=$ both $X_*$'s and $X_\bullet$'s.
Precisely, let $\mu:=\lambda_D|_{V_D\cap R_D}$. %=\bigotimes_{\tau_i\in D}\lambda_{\tau_i,s_i}|_{V_D\cap R_D}$
Then $\mu|_{X_{\tau_i}}{=}\phi_{\tau_i,s_i}$ where $\tau_i\in\{\al_{1,4},\al_{2,5},\al_{5,6},\al_{6,7},\al_{7,8}\}=*'s$
%i.e. $\mu|_{X_{\al_{1,4}}}{=}\phi_{\al_{1,4},s_1}$, $\mu|_{X_{\al_{2,5}}}{=}\phi_{\al_{2,5},s_2}$, $\mu|_{X_{\al_{5,6}}}{=}\phi_{\al_{5,6},s_3}$, $\mu|_{X_{\al_{6,7}}}{=}\phi_{\al_{6,7},s_4}$, $\mu|_{X_{\al_{7,8}}}{=}\phi_{\al_{7,8},s_5}$ for
and $s_i{\in}\F_q^\times$.
Notice that $\mu$ is a linear character of $R_D\cap V_D\subset U_9(q)$.

Since $K_D\cap R_D\subset \ker(\mu)\cap \ker({\mu}^{R_D})$, we proceed the induction of $\mu$ in the quotient group $R_D/(K_D\cap R_D)$.  Let $R:=R_D/(K_D\cap R_D)$ and $V:=(V_D\cap R_D)/(K_D\cap R_D)$. Abusing terminology we call the image of a root subgroup of $U$ under the natural projection to its factor group by a pattern subgroup, a root subgroup.

To obtain all constituents $\chi$ of ${\mu}^{R}$, we use the following strategy. We start with the factor group $R/H$ where $H$ is the largest pattern subgroup contained in $\ker(\chi)$. Let $Q\leq R/H$ such that $Q$ is maximal in the set of all pattern subgroups $P$ of $R/H$ satisfying $[P,P]\leq\ker(\mu)$.
By Lemma \ref{lem:ExtChar}, $\mu$ extends to $Q$. Let $\lambda$ be an extension of $\mu$ to $Q$. In most of the cases,
$Q\lhd R/H$ and %the inertia group
$I_{R/H}(\lambda)=Q$ except the subcase e.

\begin{proposition}
\label{decompose_lemma}
${\mu}^{R}$ decomposes into
\begin{itemize}
\item $q^3(2q-1)$ irreducibles of degree $q^3$, each has multiplicity $q^3$;

\item $q(q-1)(3q^3+q^2+q-1)$ irreducibles of degree $q^4$, each has multiplicity $q^4$;

\item $q^2(q{+}2)(q{-}1)^2$ irreducibles of degree $q^5$, each has multiplicity $q^5$.

\end{itemize}

\end{proposition}

\begin{proof}
Since $X_{\al_{1,3}} X_{\al_{3,8}}\subset Z(R)$, ${\mu}^{L_1}$ splits into $q^2$ linears %characters
of $L_1:=VX_{\al_{1,3}}X_{\al_{3,8}}$. Let $\lambda_1$ be an extension of $\mu$ to $L_1$. We divide the proof into four cases:
\begin{itemize}
\item[] Case 1: $\lambda_1(X_{\al_{1,3}})\neq \{1\} \neq \lambda_1(X_{\al_{3,8}})$, there are $(q-1)^2$ such $\lam_1$'s.

\item[] Case 2: $\lambda_1(X_{\al_{1,3}})\neq \{1\} = \lambda_1(X_{\al_{3,8}})$, there are $(q-1)$ such $\lam_1$'s.

\item[] Case 3: $\lambda_1(X_{\al_{1,3}})= \{1\} \neq \lambda_1(X_{\al_{3,8}})$,  there are $(q-1)$ such $\lam_1$'s.

\item[] Case 4: $\lambda_1(X_{\al_{1,3}})= \{1\} = \lambda_1(X_{\al_{3,8}})$, there is only one such $\lam_1$.

\end{itemize}

\textbf{Case 1}: Let $L_2:=L_1X_{\al_{1,2}}X_{\al_2}X_{\al_{2,3}}X_{\al_{2,4}}X_{\al_6}X_{\al_{3,7}}X_{\al_{4,7}}X_{\al_{4,8}}$. Check that $L_2\unlhd R$, $[L_2,L_2]\leq\ker(\lambda_1)$ and $R=L_2X_{\al_1}X_{\al_3}X_{\al_5}X_{\al_7}X_{\al_8}$. So $\lambda_1^{L_2}$ splits into $q^8$ linear characters
of $L_2$. Let $\lambda_2$ be an extension of $\lambda_1$ to $L_2$. Check that the inertia group $I_{R}(\lambda_2)=L_2$ and, by Clifford theory, $\lambda_2^{R}\in \Irr(R)$ has degree $q^5$. %$|\{{\lambda_2}^x:x\in X_{\al_1}X_{\al_3}X_{\al_6}X_{\al_8}X_{\al_{2,4}}\}|=
So we obtain $(q-1)^2q^3$ irreducible constituents of degree $q^5$, each has multiplicity $q^5$.

\begin{remark} %simplify the work
To avoid more notations, we still denote by $R$ its quotient groups $R/H$ where $H$ is some normal pattern subgroup of $R$ in the coming proofs.
\end{remark}

\textbf{Case 2}: Here $X_{\al_{3,8}}\subset \ker(\lam_1^{R})$, we work with $R/X_{\al_{3,8}}$.
Since $X_{\al_{4,8}}\subset Z(R)$, $\lambda_1$ extends to
%$q$ linear characters of
$H_2:=L_1X_{\al_{4,8}}$. Let $\eta_2$ be an extension of $\lambda_1$ to $H_2$.

If $\eta_2(X_{\al_{4,8}})\neq \{1\}$, then let $H_3:=H_2X_{\al_2}X_{\al_6}X_{\al_8}X_{\al_{1,2}}X_{\al_{2,3}}X_{\al_{2,4}}X_{\al_{3,7}}$. Check that $[H_3,H_3]\leq \ker(\eta_2)$, $H_3\unlhd R$ and $R=H_3X_{\al_1}X_{\al_3}X_{\al_5}X_{\al_7}X_{\al_{4,7}}$. So $\eta_2^{H_3}$ splits into $q^7$ linear characters. Call $\eta_3$ an extension of $\eta_2$ to $H_3$. Check that $I_{R}(\eta_3)=H_3$ and $\eta_3^{R}\in \Irr(R)$. Thus there are $(q-1)^2q^2$ irreducible constituents of degree $q^5$, each has multiplicity $q^5$.

Otherwise, if  $\eta_2(X_{\al_{4,8}})= \{1\}$, then $X_{\al_{4,8}}\subset \ker({\eta_2}^{R})$. We work with $R/X_{\al_{4,8}}$.
Let $H_3{:=}H_2X_{\al_2}X_{\al_6}X_{\al_8}X_{\al_{1,2}}X_{\al_{2,3}}X_{\al_{2,4}}X_{\al_{3,7}}X_{\al_{4,7}}$. Check that $[H_3,H_3]{\leq} \ker(\eta_2)$, $H_3\unlhd R$ and $R=H_3X_{\al_1}X_{\al_3}X_{\al_5}X_{\al_7}$. So $\eta_2^{H_3}$ splits into $q^8$ linear characters. Call $\eta_3$ an extension of $\eta_2$ to $H_3$. Check that $I_{R}(\eta_3)=H_3$ and $\eta_3^{R}\in \Irr(R)$. Thus, there are $(q-1)q^4$ irreducible constituents of degree $q^4$, each has multiplicity $q^4$.

\smallskip
\textbf{Case 3}: Here $X_{\al_{1,3}}\subset \ker(\lam_1^{R})$, we work with $R/X_{\al_{1,8}}$. Since $X_{\al_{1,2}}\subset Z(R)$, $\lam_1$ extends to $N_2:=L_1X_{\al_{1,2}}$. Call $\lambda_2$ an extension of $\lambda_1$ to $N_2$.

If $\lam_2(X_{\al_{1,2}})\neq \{1\}$, let $N_3:=N_2X_{\al_1}X_{\al_5}X_{\al_7}X_{\al_{2,3}}X_{\al_{3,7}}X_{\al_{4,7}}X_{\al_{4,8}}$. Check that $[N_3,N_3]\leq\ker(\lambda_2)$, $N_3\unlhd R$ and $R=N_3X_{\al_2}X_{\al_3}X_{\al_6}X_{\al_8}X_{\al_{2,4}}$. So $\lambda_2^{N_3}$ splits into $q^7$ linear characters. Let $\lambda_3$ be an extension of $\lambda_2$ to $N_3$. Check that $I_{R}(\lambda_3)=N_3$ and $\lambda_3^{R}\in \Irr(R)$. Thus there are $(q-1)^2q^2$ irreducible constituents of degree $q^5$, each has multiplicity $q^5$.

Otherwise, if $\lambda_2(X_{\al_{1,2}})= \{1\}$,   then $X_{\al_{1,2}}\subset \ker(\lambda_2)$. We work with $R/X_{\al_{1,2}}$. Let $N_3{:=}N_2X_{\al_1}X_{\al_2}X_{\al_5}X_{\al_7}X_{\al_{2,3}}X_{\al_{3,7}}X_{\al_{4,7}}X_{\al_{4,8}}$. Check that $[N_3,N_3]{\leq}\ker(\lambda_2)$, $R=N_3X_{\al_3}X_{\al_6}X_{\al_8}X_{\al_{2,4}}$ and $N_3\unlhd R$. So $\lambda_2^{N_3}$ splits into $q^8$ linear characters. Call $\lambda_3$ an extension of $\lambda_2$ to $N_3$. Check that $I_{R}(\lambda_3)=N_3$ and $\lambda_3^{R}\in \Irr(R)$. Thus there are $(q-1)q^4$ irreducible constituents of degree $q^4$, each has multiplicity $q^4$.

\smallskip
\textbf{Case 4}: Here $X_{\al_{1,3}} X_{\al_{3,8}}\subset \ker(\lam_1^{R})$, we work with $R/X_{\al_{1,3}} X_{\al_{3,8}}$. Since $X_{\al_{1,2}}$
$X_{\al_{2,3}} X_{\al_{3,7}} X_{\al_{4,8}}\subset Z(R)$, $\lam_1$ extends to $T_2:=L_1X_{\al_{1,2}}X_{\al_{2,3}}X_{\al_{3,7}}X_{\al_{4,8}}=Z(R)$. Let $\mu_2$ be an extension of $\lambda_1$ to $T_2$. Notice that there are nine pattern special $p$-groups $S_i$ of type $q^{1+2}$ whose centers are contained in $Z(R)$; for each $X_\al\subset R-Z(R)$, there exists exactly a pair $S_i\neq S_j$ such that  $S_i\cap S_j=X_\al$. %as follows:
\begin{center}
%$Z(R)=X_{\al_{1,2}}X_{\al_{2,3}}X_{\al_{3,7}}X_{\al_{4,8}}X_{\al_{7,8}}X_{\al_{6,7}}X_{\al_{5,6}}X_{\al_{2,5}} X_{\al_{1,4}}$ and \\
$\begin{array}{lll}
S_1=X_{\al_1}X_{\al_{1,2}} X_{\al_2}, &   S_2= X_{\al_2}X_{\al_{2,3}}X_{\al_3}, &  S_3=X_{\al_3}X_{\al_{3,7}}X_{\al_{4,7}},
\\
 S_4=X_{\al_{4,7}} X_{\al_{4,8}} X_{\al_8} ,& S_5= X_{\al_8} X_{\al_{7,8}} X_{\al_7} ,&    S_6= X_{\al_7} X_{\al_{6,7}} X_{\al_6},
\\
 S_7=X_{\al_6} X_{\al_{5,6}} X_{\al_5} ,&   S_8= X_{\al_5} X_{\al_{2,5}} X_{\al_{2,4}} ,& S_9= X_{\al_{2,4}} X_{\al_{1,4}} X_{\al_1}.
\end{array}$
\end{center}

We shall study five subcases based on which of the four root subgroups $X_{\al_{1,2}}$, $X_{\al_{2,3}}$, $X_{\al_{3,7}}$, $X_{\al_{4,8}}$ are contained in $\ker(\mu_2)$.

\smallskip
{\em Subcase a:} All are contained in $\ker(\mu_2)$. Let $T_3:=T_2X_{\al_2}X_{\al_3}X_{\al_{4,7}}X_{\al_1}X_{\al_5}X_{\al_7}$. Since $[T_3,T_3]\leq \ker(\mu_2)$, $\mu_2$ extends to $T_3$. Call $\mu_3$ an extension of $\mu_2$ to $T_3$. Each $\mu_3$ induces irreducibly to $R$ by checking $I_{R}(\mu_3)=T_3$. Thus, there are $q^3$ irreducible constituents of degree $q^3$, each has multiplicity $q^3$.

\smallskip
{\em Subcase b:} Three of them are in $\ker(\mu_2)$ and one is not. There are ${4\choose1}=4$ smaller cases corresponding to each root subgroup $\not\subset\ker(\mu_2)$. Let $\bar{R}$ denote the factor group of $R$ by the three root subgroups $\subset\ker(\mu_2)$. They fall into two families where $|Z(\bar{R})|=$ $q^8$ or $q^7$. If either $X_{\al_{1,2}}$ or $X_{\al_{4,8}}\not\subset\ker(\mu_2)$, then $|Z(\bar{R})|=q^8$ and $\mu_2^{\bar{R}}$ has $q^3(q-1)$ irreducible constituents of degree $q^3$, each has multiplicity $q^3$. Otherwise, if either $X_{\al_{2,3}}$ or $X_{\al_{3,7}}\not\subset\ker(\mu_2)$, then $|Z(\bar{R})|=q^7$ and $\mu_2^{\bar{R}}$ has $q(q-1)$ irreducible constituents of degree $q^4$, each has multiplicity $q^4$. We shall decompose $\mu_2^R$ for one case of each type, the others are similar.

If $X_{\al_{1,2}}\not\subset \ker(\mu_2)$, let $T_3:=T_2X_{\al_3}X_{\al_{4,7}}X_{\al_2}X_{\al_{2,4}}X_{\al_6}X_{\al_8}$. If $X_{\al_{2,3}}\not\subset \ker(\mu_2)$, let $T_3:=T_2X_{\al_{4,7}}X_{\al_2}X_{\al_{2,4}}X_{\al_6}X_{\al_8}$. Since $[T_3,T_3]\leq\ker(\mu_2)$, $\mu_2$ extends to $T_3$. Call $\mu_3$ an extension of $\mu_2$ to $T_3$. Then $\mu_3^R\in\Irr(R)$ by checking $I_{R}(\mu_3)=T_3$.

So we obtain $2(q-1)q^3$ irreducible constituents of degree $q^3$, each has  multiplicity $q^3$, and $2(q-1)q$ irreducible constituents of degree $q^4$, each has  multiplicity $q^4$.

\smallskip
{\em Subcase c:} Two of them are in $\ker(\mu_2)$ and the others are not. There are ${4\choose 2}=6$ smaller cases. Applied the same argument, in each case $\mu_2^R$ has $q(q-1)^2$ irreducible constituents of degree $q^4$, each has multiplicity $q^4$. Thus, we obtain $6(q-1)^2q$ irreducible constituents of degree $q^4$, each has multiplicity $q^4$.

\smallskip
{\em Subcase d:} One is in $\ker(\mu_2)$ and the other three are not. There are 4 smaller cases. WLOG, we suppose $X_{\al_{1,2}}{\subset} \ker(\mu_2)$. Let $T_3{:=}T_2X_{\al_1}X_{\al_5}X_{\al_7}X_{\al_{4,7}}X_{\al_2}$. Since $[T_3,T_3]\leq\ker(\mu_2)$, $\mu_2$ extends to $T_3$. Call $\mu_3$ an extension of $\mu_2$. Then $\mu_3^R\in\Irr(R)$ by checking $I_{R}(\mu_3)=T_3$. So there are $4(q-1)^3q$ irreducible constituents of degree $q^4$, each has multiplicity $q^4$.

\smallskip
{\em Subcase e:} None of them is in $\ker(\mu_2)$. Let $Q_3:=T_2X_{\al_1}X_{\al_5}X_{\al_7}X_{\al_{4,7}}\unlhd R$. Since $[Q_3,Q_3]\leq\ker(\mu_2)$, $\mu_2$ extends to $Q_3$. Call $\mu_3$ an extension of $\mu_2$ to $Q_3$.
Let $\mu_3|_{X_{\beta_i}}=\phi_{\beta_i,s_i}$ where $\beta_i\in\{\al_{1,2},\al_{2,3},\al_{3,7},\al_{4,8},\al_{7,8},\al_{6,7},\al_{5,6},\al_{2,5},\al_{1,4}\}$ and $s_i\in\F_q^\times$. The inertia group $I_{R}(\mu_3)$ is generated by  $Q_3$ and $x(a)$ for all $a\in\F_q$ where \[x(a):=x_{\al_2}(-\frac{s_{1,4}}{s_{1,2}}a)x_{\al_6}(\frac{s_{2,5}}{s_{5,6}}a)x_{\al_8}(\frac{s_{2,5}s_{6,7}}{s_{5,6}s_{7,8}}a)
x_{\al_3}(\frac{s_{2,5}s_{4,8}s_{6,7}}{s_{3,7}s_{5,6}s_{7,8}}a)x_{\al_{2,4}}(a).\]

We claim $|I_{R}(\mu_3):Q_3|=q$ and $[I_{R}(\mu_3),I_{R}(\mu_3)]\leq \ker(\mu_3)$ by showing that $[x_\beta(c),x(a)]\in \ker(\mu_3)$ where $\beta\in\{\al_1,\al_5,\al_7,\al_{4,7}\}$ and $c,a\in\F_q$. By the nilpotent class 2 of $R$, $[x,yz]=[x,z][x,y]^z=[x,z][x,y]$. E.g. with $\beta=\al_1$, we have
\[\begin{array}{ll} [x_{\al_1}(c),x(a)]&=[x_{\al_1}(c),x_{\al_{2,4}}(a)][x_{\al_1}(c),x_{\al_2}(-\frac{s_{1,4}}{s_{1,2}}a)]\\
&=x_{\al_{1,4}}(ac)x_{\al_{1,2}}(-\frac{s_{1,4}}{s_{1,2}}ac).
\end{array}\]

Therefore,
\[\begin{array}{ll}
\mu_3([x_{\al_1}(c),x(a)])&=\mu_3(x_{\al_{1,4}}(ac)x_{\al_{1,2}}(-\frac{s_{1,4}}{s_{1,2}}ac))\\
&=\phi_{\al_{1,4},s_{1,4}}(x_{\al_{1,4}}(ac))\phi_{\al_{1,2},s_{1,2}}(x_{\al_{1,2}}(-\frac{s_{1,4}}{s_{1,2}}ac))\\
&=\phi(s_{1,4}ac-s_{1,2}\frac{s_{1,4}}{s_{1,2}}ac)=1.
\end{array}\]

So $\mu_3$ extends to $I_{R}(\mu_3)$ and each extension $\mu_4$ induces irreducibly to $R$. Thus, there are $(q-1)^4q$ irreducible constituents of degree $q^4$, each has multiplicity $q^4$.
\end{proof}

\smallskip
\begin{remark} In the subcase e, $I_{R}(\mu_3)$ is not a pattern subgroup. By Theorem \ref{main_thm} $\xi_{D,\phi}=\lam_D^{U}$ has $|\Irr(V_DR_D,\lambda_D)|{=}(q-1)^4q$ distinct constituents of degree $q^{4+12}{=}q^{16}$.
%we get $(q-1)^4q$ irreducible constituents of $\lam_D^{V_DR_D}$. Therefore, by Theorem \ref{main_lemma}, $\xi_{D,\phi}=\lam_D^{U}$ has $(q-1)^4q$ distinct constituents of degree $q^{4+12}=q^{16}$.
Since $|D|=9$, there are $(q-1)^9$ such supercharacters.  So we have $(q-1)^{13}q$ not well-induced characters as stated in \cite{Anton1}.
For $q{=}2^f$, $x(a)^2{=}x_{\al_{2,3}}(\frac{s_{1,4}s_{2,5}s_{4,8}s_{6,7}}{s_{1,2}s_{3,7}s_{5,6}s_{7,8}}a^2){\notin} \ker(\mu_4)$  and the order ${\rm o}(x(a)){=}4$ for all  $a\in\F_q^\times$. There is $a_0\in \F_q$ such that $\mu_4(x(a_0)^2){=}-1$. Thus, $\mu_4(x(a_0))=\pm i\in \C{-}\R$, i.e. $\mu_4$ is an irrational linear character of $I_{R}(\mu_3)$. This explains why $\mu_4^{R}$  remains irrational, and so does its corresponding character in $\Irr(\lambda_D^U)$. This gives a way to construct the pair of irrational irreducible characters of degree $2^{16}$ of $U_{13}(2)$ for Isaacs-Karagueuzian's Conjecture \cite{Is-Kar}.
\end{remark}

\section*{Acknowledgement}
The content of this paper was presented at Conference on Algebraic Topology, Group Theory and Representation Theory, dedicated to 60-th birthdays of Professor Ron Solomon and Professor Bob Oliver, which was held on the Isle of Skye in 2009. The author is grateful to the organizers of the conference.
%The author would like to thank the referee for giving time to correct the scripts with many helpful comments and suggestions.

\end{document}